\documentclass{amsart}
\usepackage{latexsym,amsmath,amsfonts,amssymb,amsxtra}
\usepackage{color,latexsym,amsmath,amsthm,amsfonts,amssymb,xcolor,verbatim,lipsum,layout}

\usepackage{enumerate}

\usepackage{soul,cite,lineno,cancel}
\modulolinenumbers[2]

\scrollmode

\newtheorem{thm}{Theorem}[section]
\newtheorem{deff}[thm]{Definition}

\newtheorem{rem}[thm]{Remark}
\newtheorem{prop}[thm]{Proposition}

\newcommand{\vG}{\varGamma}
\newcommand{\vD}{\varDelta}
\def\vL{\varLambda}
\newcommand{\ve}{\varepsilon}
\newcommand{\wt}{\widetilde}

\newcommand{\vO}{\varOmega}
\newcommand{\vS}{\varSigma}
\newcommand{\ov}{\overline}
\newcommand{\mg}{\marginpar}

\def\N{{\mathbb N}}
\def\R{{\mathbb R}}

\def\mbZ{\mathbb Z}

\def\mcS{{\mathcal S}}
\def\mcL{{\mathcal L}}

\def\mcZ{{\mathcal Z}}

\def\emp{\emptyset}

\def\bs{\boldsymbol}
\def\vX{\varXi}
\def\be{\begin{equation}}

\def\ee{\end{equation}}
\def\mcD{\mathcal D}
\newcommand\bas{\begin{eqnarray*}}

\newcommand\eas{\end{eqnarray*}}
\def\ba{\begin{eqnarray}}

\def\ea{\end{eqnarray}}
\def\ba*{\begin{eqnarray*}}
\def\ea*{\end{eqnarray*}}

\def\lg{\langle}
\def\rg{\rangle}

\begin{document}

\title{Conditional Expectations in  Banach spaces with RNP}

\author[K. Musia{\l}]{Kazimierz Musia\l}
\address{50-384 Wroc{\l}aw, Pl. Grunwaldzki 2/4, Poland}
\email{kazimierz.musial@math.uni.wroc.pl}
\thanks{}


\date{\today}

\dedicatory{}


\begin{abstract}
Let $X$ be a Banach space with RNP, $(\vO,\vS,\mu)$ be a complete probability space and $\vG:\vO\to{cb(X)}$ (nonempty, closed convex and bounded subsets of $X$) be a multifunction. Assume that $\vX\subset\vS$ is  a $\sigma$-algebra and the multimeasure $M$ defined by the Pettis integral of $\vG$ be such that the restriction of $M$ to $\vX$ is of $\sigma$-finite variation. Using a lifting, I prove the existence of an Effros measurable conditional expectation of $\vG$ and present its representation in terms of quasi-selections of $\vG$. I apply then the description to martingales of Pettis integrable multifunctions  obtaining a scalarly equivalent martingale of measurable multifunctions with many martingale selections. In general the situation cannot be reduced to the separable space.
\end{abstract}

\maketitle

\mg{arXiv}

\section{Introduction.}
In general investigation of multifunctions is based on the assumption of their measurability in the sense of Effros and separability of the Banach range space $X$. In particular , while investigating the problem of existence of conditional expectation of a Pettis integrable multifunction, as a rule it is assumed that the multifunction is  Effros measurable, scalarly (quasi-) integrable and the spaces $X$ and $X^*$  are separable and have RNP. The conditional expectation of $\vG$ with respect to a sub-$\sigma$-algebra is then well described by the set of conditional expectations of selections of $\vG$. This is possible due to the Castaing representation of $\vG$. The  measurability of $\vG$ is the most essential assumption because it guarantees existence of selections. The  assumption of separability of $X^*$ is postponed for multifunctions with weakly compact values. If $X$ is non-separable, then Castaing representation in general fails and one has to use a different method.

This paper is a continuation of \cite{mu25}, where I prove that several multimeasures have Effros measurable (but without the separability assumption) Radon-Nikod\'{y}m densities. In this paper, instead of separability and RNP of the range space $X$ I assume RNP of $X$ only. Such an assumption makes it possible to prove that quite arbitrary Pettis integrable multifunction which takes as its  values closed convex and bounded subsets of a Banach space with RNP, possesses the conditional expectation with respect to any sub-$\sigma$-algebra on which its integral is of $\sigma$-finite variation. Moreover, the conditional expectation can be chosen to be Effros measurable. That fact generalizes earlier known existence results. The conditional expectation is fully described by the conditional expectations of the quasi-selections of the multifunction under consideration. Then I apply the existence theorem to martingale theory. One of the consequences is that  for each martingale of Pettis integrable multifunctions there exists a martingale of measurable multifunctions that are scalarly equivalent to the initial multifunctions and possess many martingale selections.

To have a better representation of the results I introduce an extension of the ordinary lifting on the space $L_{\infty}(\mu)$  to a selection on $L_1(\mu)$ that coincides with the initial lifting on $L_{\infty}(\mu)$. 

A presentation of the general theory of Pettis integration of vector functions can be found in \cite{mu}, \cite{mu3} and \cite{ta}. For the general theory of Pettis integrable multifunctions I recommend \cite{ckr}, \cite{eh}, \cite{mu4} and \cite{mu8}.
\section{Preliminaries.}
Throughout  $(\vO,\vS,\mu)$ is a complete
probability space and $\mcL_0(\mu)$ is the space of all $\mu$-measurable and a.e. finite real valued functions. $L_0(\mu)$ is the quotient space of $\mu$-equivalent functions. $\mcL_{\infty}(\mu)$ is the space of all $\mu$-measurable and a.e. bounded real valued functions. $L_{\infty}(\mu)$ is the ordinary quotient space of $\mcL_{\infty}(\mu)$. Similarly in case of $\mcL_1(\mu)$ and $L_1(\mu)$.

A family $W\subset{L_1(\mu)}$ (or just a family $W$ of integrable functions, not equivalence classes) is {\it uniformly integrable} if $W$ is bounded in $L_1(\mu)$ and for each $\ve>0$
there exists $\delta>0$ such that if $\mu(A)<\delta$, then $\sup_{f\in{W}}\int_A|f|\,d\mu<\ve$.

$\rho$ is an arbitrary lifting on $L_{\infty}(\mu)$ and on $(\vO,\vS,\mu)$ (formally these are two different liftings but each of them uniquely determines the other). The family of all liftings on  $L_{\infty}(\mu)$  is denoted by $\vL(\mu)$.

$X$ is a Banach space with its dual $X^*$ and the closed unit ball of $X$ is denoted by $B_X$.
If $A\subset X^*$,
then $\ov{A}^*$ is the weak$^*$-closure of $A$ and $\ov{A}$ is its norm closure. If
$A\subset{X}$ is nonempty, then we write
$|A|:\,=\sup\{\|x\|:x\in{A}\}$. $\ov\int$ and $\underline\int$ denote respectively the upper and lower integrals.

$c(X)$ denotes the collection of all
nonempty closed convex subsets of $X$, $cb(X)$ is the collection
of all  bounded elements of $c(X)$, $cwk(X)$  denotes the family of all  weakly compact
elements of $cb(X)$ and  $cw^*k(X^*)$  denotes the family of all non-empty convex and weak$^*$-compact subsets of $X^*$.
 For every $C \in c(X)$ the
{\it support function of} $C$ is denoted by $s( \cdot, C)$ and
defined on $X^*$ by $s(x^*, C) = \sup \{ \langle x^*,x \rangle : \
x \in C\}$, for each $x^* \in X^*$. Notice that we have always $s(x^*, C)\neq-\infty$. $d_H$ stands for the Hausdorff distance on the hyperspace $cb(X)$.

A function $f:\vO\to\R$ is \textit{quasi-integrable} (cf. \cite{ne}), if the integral $\int_{\vO}f\,d\mu$ exists.
A map $\vG:\vO\to c(X)$ is called a {\it multifunction}.  $\vG$ is said to be {\it
scalarly measurable} ({\it
scalarly integrable, scalarly quasi-integrable}) if for every $ x^* \in X^*$, the map
$s(x^*,\vG(\cdot))$ is measurable (integrable, quasi-integrable).

A multifunction $\vG:\vO\to{c(X^*)}$ is said to be  {\it
weak$^*$-scalarly measurable} ( {\it
weak$^*$-scalarly integrable, weak$^*$-scalarly quasi-integrable})
if for every $ x\in X$, the map
$s(x,\vG(\cdot))$ is measurable (integrable, quasi-integrable).

I associate with each $\vG$ the trace of  $\vG$ on $\mcL_0(\mu)$:
$$
{\mathcal Z}_{\vG}:=\{s(x^*,\vG):\|x^*\|\leq 1\}\,,
$$
The set of equivalence classes of functions scalarly equivalent will be denoted by $\mbZ_{\vG}$.

Multifunctions $\vD,\vG:\vO\to{c(X)}[c(X^*)]$ are scalarly (weak$^*$ scalarly) equivalent, if for each $x^*\in{X^*}\;[x\in{X}]$ we have $s(x^*,\vD)=s(x^*,\vG)\;[s(x,\vD)=s(x,\vG)]$ a.e. We write then $\vD\stackrel{s}{\thickapprox}\vG\;[\vD\stackrel{w^*s}{\thickapprox}\vG]$.

Similarly, we write $\vD\subset\vG$, if $\vD(\omega)\subset\vG(\omega)$ for every $\omega\in\vO$.

A  multifunction $\vG:\vO\to c(X),\;[c(X^*)]$ is
 [weak$^*$-]scalarly bounded if there is a constant $M>0$ such that
for every $y\in X^*\,[X]$
$$
|s(y,\vG)|\leq M\|y\| \qquad a.e.
$$
Following the generally accepted terminology $\vG:\vO\to{c(X)}$ is called measurable   if
\be\label{E0}
\{\omega\in\vO\colon U\cap\vG(\omega)\neq\emp\}\in\vS\,.
\ee
 for every norm open set $U\subset{X}$.

A function $\gamma:\vO\to X$
is called a {\it selection of} $\vG:\vO\to{c(X)}$, if for almost every $\omega\in
\vO$, one has $\gamma(\omega) \in \vG(\omega)$. $\gamma:\vO\to X$ ($\gamma:\vO\to X^{**}$)
is called a {\it quasi-selection (see \cite{mu4}) of} $\vG:\vO\to{c(X)}$ (a {\it weak$^*$ quasi-selection of} $\vG$), if $\gamma$ is scalarly
measurable (weak$^*$ scalarly measurable) and for each $x^*\in{X^*}$ one has $\langle{x^*,\gamma(\omega)}\rangle\in \langle{x^*,\vG(\omega)}\rangle$ for $\mu$-almost
every $\omega\in\vO$ (the exceptional sets depend on $x^*$).  The
collection of all quasi-selections of $\vG$ (weak$^*$ quasi-selections of $\vG$) will be denoted by ${\mathcal
{QS}}_{\vG}$ (${\mathcal{W^{**}QS}}_{\vG}$). A scalarly
measurable $\gamma:\vO\to{X}$ (a weak$^*$ scalarly measurable $\gamma:\vO\to{X^{**}}$) is a (weak$^*$) quasi-selection of $\vG:\vO\to{c(X)}$ if and only if for each $x^*\in{X^*}$
\begin{equation}
\langle{x^*,\gamma}\rangle\leq{s(x^*,\vG)} \quad a.e.
\end{equation}

 A function $\gamma:\vO\to X^*$ is said to be a
weak$^*$-quasi-selection of a multifunction $\vG:\vO\to c(X^*)$ if $\gamma$ is weak$^*$-scalarly
measurable and we have $\langle{x,\gamma(\omega)}\rangle\in \langle{x,\vG(\omega)}\rangle$ for $\mu$-almost
every $\omega\in\vO$ (the exceptional sets depend on $x$).  The
collection of all weak$^*$-quasi-selections of $\vG$ will be denoted by ${\mathcal
{W^*QS}}_{\vG}$. One can easily check that a  weak$^*$-scalarly measurable
$\gamma:\vO\to{X^*}$ is  a  weak$^*$-quasi-selection of $\vG$ if and only if for each $x\in{X}$
\begin{equation}\label{E1}
\langle{x,\gamma}\rangle\leq{s(x,\vG)} \quad a.e.
\end{equation}
If $\gamma\in{\mathcal{W^{**}QS}}_{\vG}$ ($\gamma\in{\mathcal{W^*QS}}_{\vG}\,,\;\gamma\in{\mathcal{QS}}_{\vG}$) is Pettis integrable, then we will write it as $\gamma\in{\mathcal{W^{**}QS}}_{\vG}^P$ ($\gamma\in{\mathcal{W^*QS}}_{\vG}^P\,,\;\gamma\in{\mathcal{QS}}_{\vG}^P$).

If $\gamma:\vO\to{X^*}$ is a weak$^*$-scalarly measurable and weak$^*$-scalarly bounded function, then $\rho(\gamma):\vO\to{X^*}$ is defined by $\langle{x,\rho(\gamma)(\omega)}\rangle:=\rho(\langle{x,\gamma}\rangle)(\omega)$ (see \cite[ch. IV.5]{IT}). If $\gamma:\vO\to{X}$ is scalarly bounded and scalarly measurable, then $\rho(\gamma):\vO\to{X^{**}}$ is defined by $\langle{x^*,\rho(\gamma)(\omega)}\rangle:=\rho(\langle{x^*,\gamma}\rangle)(\omega)$, $x^*\in{X^*}$.

A map $M:\vS\to c(X)$ is additive, if $M(A\cup{B})=M(A)\oplus{M(B)}:=\ov{M(A)+M(B)}$ for each pair of disjoint elements of $\vS$.\\
An additive map $M:\vS\to c(X)$ is called a {\it multimeasure} if
$s(x^*,M(\cdot)):\vS\to(-\infty,+\infty]$ is a $\sigma$-finite measure, for every $x^*\in{X^*}$.
$M:\vS\to c(X^*)$ is called a {\it weak$^*$-multimeasure} if
$s(x,M(\cdot)):\vS\to(-\infty,+\infty]$ is a  $\sigma$-finite measure, for every $x\in{X}$. $\bigcup{M(\vS)}:=\bigcup_{E\in\vS}M(E)$ denotes the range of $M$ in $X$ or $X^*$. One can easily see that if $M$ is a weak$^*$-multimeasure (multimeasure), and if for $y\in{X}(X^*)$ the equality $s(y,M)=s(y,M)^+-s(y,M)^-$ is the Jordan decomposition of $s(y,M)$, then
\be
s(y,M(\vO))^-<\infty.
\ee
 If $M:\vS\to{cb(X)}$ is countably additive in the metric $d_H$, then it is called an $h$-multimeasure (cf. \cite{hp}).
 A helpful tool to study the $h$-multimeasures is the R{\aa}dstr\"{o}m embedding
 $j:cb(X)\to \ell_{\infty}(B_{X^*})$, defined
 by $j(A):=s(\cdot, A)$, (see, for example  \cite[Theorem 3.2.9 and Theorem 3.2.4(1)]{Beer} or \cite[Theorem II-19]{cv})
  It is known that $B_{X^*}$ can be embedded into $\ell_{\infty}^*(B_{X^*})$  by the mapping $x^*\longrightarrow e_{x^*}$, where
$\langle{e_{x^*},h}\rangle=h(x^*)$, for each $h\in\ell_{\infty}(B_{X^*})$. Moreover, the range of $B_{X^*}$ is a norming subset of
 $\ell_{\infty}^*(B_{X^*})$.
The embedding $j$  satisfies the following properties:
\begin{description}
\item[a)] $j(\alpha A \,  \oplus \,  \beta C) = \alpha j(A) + \beta j(C)$ for every $A,C\in  cb(X),\,\, \alpha, \beta \in  \mathbb{R}{}^+$;
\item[b)] $d_H(A,C)=\|j(A)-j(C)\|_{\infty},\quad A,C\in  cb(X)$;
\item[c)] $j(cb(X))$  is a closed cone in the space
$\ell_{\infty}(B_{X^*})$ equipped with the norm of the uniform convergence.
\end{description}
If $M$ is a point map, then we talk about measure and weak$^*$-measure,
respectively.  \\
 If $m:\vS\to{X^*}$ is a weak$^*$-measure such that
$m(A)\in{M(A)}$, for every $A\in\vS$, then $m$ is called a {\it
weak$^*$-selection} of $M$. ${\mathcal {W^*S}}(M)$ will denote the
set of all weak$^*$-countably additive selections of $M$.
Similarly, a vector measure $m:\vS\to{X}$ such that
$m(A)\in{M(A)}$, for every $A\in\vS$, is called a {\it selection}
of $M$. ${\mathcal S}(M)$ will denote the set of all countably
additive selections of $M$.

A multimeasure $M: \vS \to c(X)$ is called
	\textit{  rich } if $M(A) = \overline{\{m(A), m \in S(M) \}}$.
	By a result of Cost\'{e}, quoted in \cite[Theorem 7.9]{hess}, this is verified for  multimeasures which take as their values  weakly compact convex subsets of a Banach space
	 $X$, and $cb(X)$-valued multimeasures, when $X$ is a Banach space possessing  RNP (\cite[Th\'{e}or\`{e}me 1]{co}).

If $M$ is a weak$^*$-multimeasure or a multimeasure, then $M$ is called absolutely continuous with respect to $\mu$ (we write then $M\ll\mu$), if $\mu(E)=0$ yields  $M(E)=\{0\}$. $M$ is dominated by $\mu$ if there exists $\alpha>0$ such that $|M(E)|\leq\alpha\mu(E)$, for every $E\in\vS$.\\
If $M$ is a weak$^*$-multimeasure or a multimeasure, then the variation of $M$ is defined by
$$
\bs{|}M\bs{|}(E):\,=\sup_{\mathcal P}\left\{\sum_{A\in\mathcal P}|M(A)|:\mathcal P\;
{\rm is\; a\; finite \;measurable\; partition \; of}\;
E\right\}.
$$
One can easily see that $\bs{|}M\bs{|}$ is a measure (cf. \cite[page 863]{mu8}) and,  $M\ll\mu$ if and only if $\bs{|}M\bs{|}\ll\mu$ (that is $\mu(E)=0\Rightarrow\bs{|}M\bs{|}(E)=0$).  If $M\ll\mu$ is of $\sigma$-finite variation, then there exists a partition $(\vO_n)_n$ of $\vO$ such that $|M(E)|\leq{n}\mu(E)$ (equivalently $\bs{|}M\bs{|}(E)\leq{n\mu(E)}$) for every $E\in\vO_n\cap\vS,\,n\in\N$.
\begin{deff} \rm Denote by $\mathcal C$ a non-empty subfamily of $c(X^*)$.
A weak$^*$-scalarly quasi-integrable multifunction $\vG:\vO\to c(X^*)$ is {\it Gelfand integrable} in $\mathcal C$,  if
 for each  set $A \in \vS$ there
exists a set $M_{\vG}^G(A)\in {\mathcal C}$  such that
\begin{equation}
s(x,M_{\vG}^G(A))=\int_A s(x,\vG)\,d\mu\,,\;\mbox{for every}\; x\in X.
\end{equation}
$M_{\vG}^G$ is a weak$^*$-multimeasure and $M_{\vG}^G(A)$ is called the Gelfand integral of $\vG$ on $A$: $(G)\int_E\vG\,d\mu:=M_{\vG}^G(A)$. One should notice that the Gelfand integral is uniquely determined in the family $cw^*k(X^*)$ but not in $c(X^*)$ (see \cite[Remark 3.12]{mu17}).
\end{deff}
\begin{deff} {\rm Denote by $\mathcal D$ a non-empty subfamily of $c(X)$.  A scalarly quasi-integrable multifunction $\vG\colon \vO\to {c(X)}$ is {\it Pettis
integrable}  in
 $\mathcal D$,  if
 for each $A\in\vS$ there exists a set
 $M_{\vG}(A)\in \mathcal D$ such that
\begin{equation}
s(x^*,M_{\vG}(A))=\int_A s(x^*,\vG)\,d\mu\quad {\rm for\;every\;}
x^*\in{X^*}.
\end{equation}
$M_{\vG}:\vS\to{\mathcal D}$ is a multimeasure and $M_{\vG}(A)$ called the {\it Pettis integral of} $\vG$ over
$A$. We write  $(P)\!\!\int_A\vG\,d\mu:=M_{\vG}(A)$. ${\mathbb P}(\mu,\mcD)$ denotes the space of all multifunctions $\vG\colon \vO\to {c(X)}$ that are Pettis integrable in $\mcD$, with convention that scalarly equivalent multifunctions are identified.
${\mathbb P}(\mu,cb(X))$ is endowed with the Pettis metric
\[
d_P(\vG,\vD):=\sup_{x^*\in{B_{X^*}}}\int_{\vO}|s(x^*,\vG)-s(x^*,\vD)|\,d\mu
\]
generating the Pettis norm
\[
\|\vG\|_P:=\sup_{x^*\in{B_{X^*}}}\int_{\vO}|s(x^*,\vG)|\,d\mu.
\]
We have (see \cite[Lemme 7]{cv})
\begin{equation}\label{E3}
\sup_{E\in\vS}d_H(M_{\vG}(E),M_{\vD}(E))\leq d_P(\vG,\vD)\leq 2\,\sup_{E\in\vS}d_H(M_{\vG}(E),M_{\vD}(E)).
\end{equation}
If $M_{\vG}$ is countably additive in the Hausdorff metric (i.e. it is an $h$-multimeasure), then we say that $\vG$ is {\it strongly Pettis integrable}.}
\end{deff}

The following well known property of lifting will be applied several times in this paper:
\begin{prop}\label{p4} Let $\rho$ be a lifting on the complete $(\vO,\vS,\mu)$ and  $\{A_t:t\in{T}\}\subset \vS$ be a collection of sets such that $A_t\subset\rho(A_t)$, for each $t\in{T}$. Then  $\bigcup_{t\in{T}}A_t\in\vS$ and $\bigcup_{t\in{T}}A_t\subset\rho(\bigcup_{t\in{T}}A_t)$.

If $\{f_t:t\in{T}\}$ is a uniformly bounded family of real-valued functions such that $f_t\leq\rho(f_t)$  for each $t\in{T}$, then $\sup_{t\in{T}}f_t$  is a measurable function and $\sup_{t\in{T}}f_t\leq\rho(\sup_{t\in{T}}f_t)$ everywhere.
\end{prop}
\section{Pseudo-lifting on $L_0(\mu)$.}
In order to have a better shape of  representations of multimeasures as Pettis integrals I introduce an extension of lifting from $L_{\infty}(\mu)$ to the whole $L_0(\mu)$. It is well known that no such an extension has the same properties as the ordinary lifting. This is in fact a selection $\rho^\#:L_0(\mu)\to\mcL_0(\mu)$ uniquely determined by the lifting $\rho$ on $L_{\infty}(\mu)$. I introduce it not because of excelent properties but because it makes possible to write down several formulae in a consistent form. If $f\in\mcL_0(\mu)$, then the set $Dom(f):=\bigcup\{\rho(E):f\;\mbox{is a.e. bounded on }E\}=\bigcup_n\rho(E_n)$ with $E_n:=\{\omega:\vert{f(\omega)}\vert\leq {n}\}$, will be called the domain of $f$.
\begin{deff} \rm
If $\rho$ is a lifting on $L_{\infty}(\mu)$, then I call the mapping $\rho^\#:\mcL_0(\mu)\to{{\mcL}_0(\mu)}$ (or $\rho\#:L_0(\mu)\to{{\mcL}_0(\mu)}$) a pseudo-lifting generated by $\rho$, if the following conditions are fulfilled:
\begin{enumerate}
\item
$\rho^\#(f):Dom(f)\to\R$ for every $f\in {{\mcL}_0(\mu)}$;
\item
$
\rho^\#(f)\chi_{\rho(E)}=\rho(f\chi_E)
$
for every $E\in\vS$ such that $f$ is a.e. bounded on $E$.
\end{enumerate}
\end{deff}
\begin{rem}\label{r2} \rm
 If $f=g\;\mu$-a.e, then $\rho^\#(f)=\rho^\#(g)$ and so  $\rho^\#$ is well defined.

If $f$ and $g$ are a.e. bounded on $E\in\vS$ and $\alpha,\beta\in\R$, then $\rho^\#(\alpha{f}+\beta{g})=\alpha\rho^\#(f)+\beta\rho^\#(g)$ on $\rho(E)$ but in general
 $\rho^\#$ is not additive. One can say only that $\rho^\#(f+g)=\rho^\#(f)+\rho^\#(g)$ a.e. More precisely, the pointwise equality holds true on $Dom(f)\cap{Dom(g)}$. If $f$ is a.e. bounded on $E\in\vS$, then
\[
\rho^\#(f)=\rho(f\chi_E)+\rho^\#(f\chi_{E^c}).
\]
If $E\in\vS$ is arbitrary, then
\[
\rho^\#(f)=\rho^\#(f\chi_E)+\rho^\#(f\chi_{E^c}).
\]
\end{rem}
\begin{deff}\rm
Let $\vG:\vO\to{cb(X^*)}$ (in particular $\gamma:\vO\to{X^*}$) be weak$^*$ scalarly measurable and $Dom(\vG):=\bigcup\{\rho(E):\vG\;\mbox{is weak$^*$ scalarly bounded on }E\}=\bigcup_n\rho(F_n)$, where $\vG$ is weak$^*$ scalarly bounded by $n$ on $F_n$ and $F_n$ is maximal in the following sense: $\mu(F_n)=\sup\{\mu(E): \forall\;x\in{B_X}\;|s(x,\vG)|\leq{n}\;\mu-a.e.\;\mbox{on }F_n\}$.  According to \cite[Proposition 1.2]{mu4} $\mu(Dom(\vG)=1$.  We define a $cb(X^*)$-valued $\rho^\#\vG=\rho^\#(\vG):Dom(\vG)\to{cb(X^*)}$ by setting
\[
\rho^\#(\vG)\chi_{\rho(E)}:=\rho(\vG\chi_E)\quad \mbox{if }\vG\;\mbox{is weak$^*$ scalarly bounded on }E\,,
\]
where $\rho(\vG\chi_E)$ is defined by \cite[Proposition 4.1]{mu17}.
\end{deff}
\begin{rem}\label{r4} \rm
 If $\vG$ is  weak$^*$ scalarly bounded on $E$, then
\begin{equation}\label{E2}
\rho^\#\vG(\omega)=\rho(\vG\chi_E)(\omega)+\rho^\#(\vG\chi_{E^c})(\omega)\quad \mbox{if } \omega\in{Dom}(\vG) .
\end{equation}
If $\vG,\vD:\vO\to{cb(X^*)}$  are  weak$^*$ scalarly bounded on $E\in\vS$ and $a,b\in[0,+\infty)$, then
$\rho^\#(a\vG\oplus {b\vD})=a\rho^\#\vG\oplus{b}\rho^\#\vD$ on $\rho(E)$.

In general one cannot say that $s(x,\rho^\#\vG)=\rho^\#[s(x,\vG)]$ everywhere.  But we have the equality
\[
s(x,\rho^\#\vG)=\rho^\#[s(x,\vG)]\quad\mbox{everywhere on }Dom(\vG)\,.
\]
In particular, if $\gamma\in\mathcal{W^*QS}_{\vG}$, then $\rho^\#(\gamma)\in\mathcal{W^*S}_{\rho^\#\vG}$ and
\begin{equation}\label{E51}
\lg{x,\rho^\#(\gamma)}\rg=\rho^\#\lg{x,\gamma}\rg \quad\mbox{everywhere on }Dom(\vG)\,.
\end{equation}
If $\gamma,\delta\in\mathcal{W^*QS}_{\vG}$ and $a,b\in\R$, then $\rho^\#(\gamma),\rho^\#(\delta)\in\mathcal{W^*S}_{\rho^\#\vG}$ and
$\rho^\#(a\gamma+b\delta)=a\rho^\#(\gamma)+b\rho^\#(\delta)$ on $Dom(\gamma)\cap{Dom(\delta)}\supset{Dom(\vG)}$.
\end{rem}
\section{Radon-Nikod\'{y}m Theorem.}
The subsequent result is a generalization of \cite[Proposition 4.3]{mu17}.
\begin{prop}\label{p14}
 If $\vG:\vO\to{cb(X^*)}$ is a weak$^*$-scalarly integrable multifunction and $M:\vS\to{cb(X^*)}$ is its Gelfand integral, then the  multifunctions $\rho^\#\vG:\vO\to{cw^*k(X^*)}$ and $\vG_{\rho^\#}:\vO\to{cb(X^*)}$ are defined for every $\omega\in{Dom(\vG)}$ by the formulae
\begin{equation}
\vG_{\rho^\#}(\omega):=\ov{conv} \left\{\rho^\#(\gamma)(\omega):
\gamma\in{\mathcal {W^*QS}}_{\vG}\right\}= \ov{\left\{\rho^\#(\gamma)(\omega):
\gamma\in{\mathcal {W^*QS}}_{\vG}\right\}}
\end{equation}
\begin{equation}\label{E5}
(\vG_M)_{\rho^\#}(\omega):=\ov{\left\{\rho^\#(\gamma)(\omega):(G)\!\!\int_{\bullet}\rho^\#(\gamma)\,d\mu\in{\mathcal{W^*S}}(M)\right\}}.
\end{equation}
Moreover, $\rho^\#\vG,\;\rho^\#(\vG_M):\vO\to{cw^*k(X^*)}\cup\emp$ are defined by
\begin{equation}
\rho^\#\vG(\omega):=\ov{\vG_{\rho^\#}(\omega)}^{\,*}
\end{equation}
\begin{equation}\label{E6}
 \rho^\#(\vG_M)(\omega)=\ov{(\vG_M)_{\rho^\#}(\omega)}^{\,*}
\end{equation}
$\vG_{\rho^\#},\;(\vG_{M^*})_{\rho^\#},\;\rho^\#\vG$ and $\rho^\#(\vG_{M^*})$
are  Gelfand integrable densities of $M$ such that
\begin{equation}\label{E7}
s(x,\vG_{\rho^\#}(\omega))=s(x,\rho^\#\vG(\omega))=\rho^\#(s(x,\vG))(\omega)\quad\mbox{for every} \;x\in{X}\;\mbox{and every}\;\omega\in{Dom(\vG)}.
\end{equation}
\end{prop}
\begin{proof}
If $\vO=\bigcup_nE_n$ and $\vG$ is weak$^*$-scalarly bounded by $n$ on $E_n$, then on each $E_n$  the Proposition coincides with \cite{mu17}[Proposition 4.3]. The rest of the proof is standard.
\end{proof}
\begin{rem}\label{r1}\rm It follows from \eqref{E7} that ${\mathcal {W^*QS}}_{\vG}={\mathcal {W^*QS}}_{\rho^\#\vG}={\mathcal {W^*QS}}_{\vG_{\rho^\#}}$.
\end{rem}
The situation becomes more complicated when we investigate $cb(X)$-valued multifunctions and $X$ is not necessarily a conjugate Banach space. We will use a result from \cite{ed} (see also \cite[Lemma 3.2]{mu25}): If  $f:\vO\to{X}$ is a scalarly bounded strongly measurable function, then $\rho(f)$ is a strongly measurable function with values in $X^{**}$ but  almost all values of $\rho(f)$ are in $X$.

If a multimeasure $M$ is defined on $\vS$, then we denote by $\mcD_M^P$ the collection of all  Pettis integrable densities of selections of $M$ and by $\mcD_M^{mP}$ the collection of all strongly measurable  members of $\mcD_M^P$.
\begin{deff}
Let $M:\vO\to{cb(X)}$ be a multimeasure dominated by $\mu$. If $\mcD_M^P\neq\emp$, then we define $\vG_{M\rho}:\vO\to{cb(X^{**})}$ by
\[
\vG_{M\rho}(\omega):=\ov{conv}\left\{\rho(\gamma)(\omega)\colon \gamma\in\mcD_M^P\right\}
=\ov{\left\{\rho(\gamma)(\omega): \gamma\in\mcD_M^P\right\}}\subset (\vG_M)_{\rho}(\omega)
\]
and  $\vG_{M\rho}^X:\vO\to{cb(X)}\cup\emp$ by
\[
\vG_{M\rho}^X(\omega):=\ov{\left\{\rho(\gamma)(\omega): \gamma\in\mcD_M^P\right\}\cap{X}}
\subset \vG_{M\rho}(\omega)\cap{X}\,.
\]
\end{deff}
It is clear that $\mu\{\omega\in\vO:\vG_{M\rho}^X(\omega)=\emp\}=0$.
If $X$ has RNP, then each element of $\mcD_M^P$ is scalarly equivalent to an element of $\mcD_M^{mP}$ and so everywhere in the above definition $\mcD_M^P$  may be replaced by $\mcD_M^{mP}$.
\begin{prop}\label{p7} \cite[Proposition 3.7]{mu25}
If $M:\vS\to{cb(X)}$ is a multimeasure dominated by $\mu$, then  $\vG_{M\rho}$ is weak$^*$ scalarly measurable and for each $x^*\in{X^*}$
\begin{equation}\label{E8}
s(x^*,\vG_{M\rho})\leq\rho[s(x^*,\vG_{M\rho})]
\end{equation}
everywhere. Moreover, the function $|\vG_{M\rho}|$ is measurable.

If  $X$ has RNP and $\rho(\vG_{M^{**}})\,,(\vG_{M^{**}})_{\rho}$ are  defined by \eqref{E5} and \eqref{E6},  then the multifunctions

$\vG_{M\rho}^X\,,\vG_{M\rho}\,,\rho(\vG_{M^{**}})\,,(\vG_{M^{**}})_{\rho}$ are weak$^*$ scalarly equivalent and
\[
s(x^*,\vG_{M\rho}^X)\leq s(x^*,\vG_{M\rho})\leq\rho[s(x^*,\vG_{M\rho}^X)]=\rho[s(x^*,\vG_{M\rho})]\notag
\]
everywhere on the set $\{\omega\in\vO: \vG_{M\rho}^X(\omega)\neq\emp\}$. Moreover, the function $|\vG_{M\rho}^X|$ is measurable.

If moreover $M$ is the Pettis integral of $\vG:\vO\to{cb(X)}$, then $\vG$ is scalarly equivalent to $\vG_{M\rho}^X$ and $\vG_{M\rho}\cap{X}$.
\end{prop}

The subsequent result is an extension of \cite[Proposition 4.8]{mu17} but under additional assumptions and is a generalization of \cite[Theorem 4.1]{ckr3}.
\begin{thm}\label{t6} \cite[Theorem3.8]{mu25}
Let $M:\vS\to{cb(X)}$ be a $\mu$-continuous multimeasure of $\sigma$-finite variation. If  $X$ has RNP and $\rho\in\vL(\mu)$, then there exists a measurable multifunction $\vG:\vO\to{cb(X)}$ such that
\begin{equation}\label{E9}
M(E)=\ov{\biggl\{(P)\!\!\int_E\gamma\,d\mu\colon \gamma\in\mcD_M^{mP}\biggr\}}=(P)\!\!\int_E\vG\,d\mu\quad\mbox{for every }E\in\vS
\end{equation}
and
\begin{equation}\label{E10}
\{\omega\in\vO\colon \vG(\omega)\cap{U}\neq\emp\}\subset\rho\{\omega\in\vO\colon \vG(\omega)\cap{U}\neq\emp\}
\end{equation}
for every norm open $U\subset{X}$.

If $M$ is dominated by $\mu$, then one may set  $\vG=X\cap\vG_{M\rho}$ or $\vG=\vG_{M\rho}^X$ and in each case the inequality $s(x^*,\vG)\leq\rho(s(x^*,\vG))$ holds true everywhere on the set $\{\omega:\vG_{M\rho}^X(\omega)\neq\emp\}$  and for every $x^*\in{X^*}$.
\end{thm}
\section{Conditional expectation.}
The problem of existence of conditional expectation of  Pettis integrable multifunctions and convergence of Pettis integrable martingales of multifunctions have been investigated in several papers. However the main results were obtained under the assumption of separability and RNP of $X$, separability of $X^*$ and Castaing representation of measurable multifunctions. If $X$ is non-separable, then there are several measurable multifunctions for which the classical Castaing representation fails. In this and the next chapter we assume RNP of $X$  but nothing is assumed about $X^*$. Under these assumptions we obtain an uncountable Castaing type representation of measurable multifunctions.

The following result is essential for our description of the conditional expectation of a multifunction:
\begin{prop}\label{p3}
Let $X$ be an arbitrary Banach space and $M:\vS\to{cb(X)}$ be a multimeasure. Moreover, let $\Xi\subset\vS$ be a $\sigma$-algebra and $M_1:=M|_{\Xi}$. If $M$ is rich, then $M_1$ is also rich  and for each $m_1\in\mcS(M_1)$ there exists an extension of $m_1$ to the whole $\vS$ being a selection of $M$.
\end{prop}
\begin{proof}
The riches of $M_1$ is a consequence of \cite[8. Proposition 4.35]{hp} and the existence of extensions of each selection of $M_1$ to a selection of $M$ follows from the proof of \cite[Proposition 4.35]{hp}. If $m_1\in\mcS(M_1)$ and $m\in\mcS(M)$ is its $\mu$-continuous extension with $m(E)=(P)\!\!\int_E f\,d\mu$, for every $E\in\vS$, then $m_1(F)=(P)\!\!\int_F\mathbb{E}_{\vX}(f)\,d\mu$ for every $F\in\vX$.
\end{proof}
\begin{deff}\rm Let $\vG:\vO\to{c(X)}$ be a multifunction that is Pettis integrable in $c(X)$ and let $\vX\subset\vS$ be a $\sigma$-algebra that is complete with respect to $\mu|_{\vX}$. A multifunction $\mathbb{E}_{\vX}(\vG):\vO\to{c(X)}$ that is scalarly measurable with respect to $\vX$ and Pettis integrable on $\vX$ in $c(X)$ is called the conditional expectation of $\vG$ with respect to $\vX$, if
\begin{equation}\label{E14}
(P)\!\!\int_E\vG\,d\mu=(P)\!\!\int_E\mathbb{E}_{\vX}(\vG)\,d\mu\quad\mbox{for every }E\in\vX\,.
\end{equation}
Each multifunction $\vD$ that is scalarly measurable with respect to $\vX$ and satisfies the equality \eqref{E14} when placed instead of $\mathbb{E}_{\vX}(\vG)$ is a version of the conditional expectation of $\vG$ with respect tom $\vX$. I will write sometimes that $\vD\in{\mathbb E}_{\vX}(\vG)$.  All versions are scalarly equivalent but not necessarily almost everywhere equal.
\end{deff}
\begin{thm}\label{t5}
Assume that $\vG:\vO\to{c(X)}$ is Pettis integrable in $cb(X)$ and  $X$ has RNP.
 Assume also that $\vX$ is a sub-$\sigma$-algebra of $\vS$ that is complete with respect to $\mu|_{\vX}$ and the multimeasure $M_1:\vX\to{cb(X)}$ being the restriction of the multimeasure $M$ to $\vX$ is dominated by $\mu$ restricted to $\vX$.  If  $\rho_1$ is a lifting on $\Xi$, then each of the following formulae defines a scalarly bounded and measurable version of the conditional expectation  $\mathbb{E}_{\vX}(\vG):\vO\to{cb(X)}$ of $\vG$ with respect to $\vX$:
\begin{align}
\mathbb{E}_{\vX}(\vG)(\omega)&=
\ov{ \left\{\rho_1[\mathbb{E}_{\vX}(\gamma)](\omega)\colon \gamma\in\mcD_M^{mP}\right\}\cap{X}}\quad\mbox{for every }\omega\in\vO\label{E15}\\
\mathbb{E}_{\vX}(\vG)(\omega)&=
X\cap\ov{ \left\{\rho_1[\mathbb{E}_{\vX}(\gamma)](\omega)\colon \gamma\in\mcD_M^{mP}\right\}}\quad\mbox{for every }\omega\in\vO\label{E16}
\end{align}
Moreover, if $x^*\in{X^*}$, then
\[
s(x^*,\mathbb{E}_{\vX}(\vG)(\omega))\leq \rho_1[s(x^*,\mathbb{E}_{\vX}(\vG))](\omega)\quad\mbox{for every }\omega\in\vO
\]
and if $U\subset{X}$ is open, then
\[
\{\omega\in\vO:\mathbb{E}_{\vX}(\vG)(\omega)\cap{U}\neq\emp\}\subset \rho_1\{\omega\in\vO:\mathbb{E}_{\vX}(\vG)(\omega)\cap{U}\neq\emp\}\,.
\]
\end{thm}
\begin{proof}
By Theorem \ref{t6} there exists a measurable scalarly bounded multifunction  $\vD:\vO\to{cb(X)}$ such that $M_1(E)=(P)\int_E\vD\,d\mu$ for every $E\in\vX$, $\vD$ is  measurable with respect to $\vX$ and $\{\omega\in\vO\colon \vD(\omega)\cap{U}\neq\emp\}\subset \rho_1\{\omega\in\vO\colon \vD(\omega)\cap{U}\neq\emp\}$ for each open $U\subset{X}$.  According to Theorem \ref{t6} we may assume that $\vD=X\cap\vD_{M_1\rho_1}$ or $\vD=\vD_{M_1\rho_1}^X$.
That proves the first equalities in \eqref{E15}-\eqref{E16}.

We have to prove  validity of \eqref{E16} and \eqref{E15}. We will concentrate on \eqref{E16}. According to Theorem \ref{t6} and Proposition \ref{p3} a measurable version of $\mathbb{E}_{\vX}(\vG)$ can be described by the formula
\begin{equation}
\mathbb{E}_{\vX}(\vG)(\omega)=X\cap\vD_{M_1\rho_1}(\omega)=X\cap\ov{\{ \rho_1(\xi)(\omega)\colon \xi\in\mcD_{M_1}^{mP}\}}.
\end{equation}

Assume now that $\xi\in \mcD_{M_1}^{mP}$ and let
\[
\varTheta(\omega):=\ov{ \left\{\rho_1[\mathbb{E}_{\vX}(\gamma)](\omega)\colon \gamma\in\mcD_M^{mP}\right\}}\,.
\]
Let $\nu_{\xi}:\vX\to{X}$ be the vector measure defined by $\nu_{\xi}(E)=(P)\!\!\int_E\xi\,d\mu=(P)\!\!\int_E\rho_1(\xi)\,d\mu$ for every $E\in\vX$.
Due to  Proposition \ref{p3} there exists  $\wt\nu_{\xi}\in\mcS(M)$ such that $\nu_{\xi}=\wt\nu_{\xi}|_{\Xi}$. Since $X$ has RNP, there exists a strongly measurable $\gamma$ such that
$\wt\nu_{\xi}(A)=(P)\!\!\int_A\gamma\,d\mu$, for every $A\in\vS$. Due to strong measurability of $\xi$ and $\gamma$ we may assume that ${\mathbb{E}}_{\vX}(\gamma)=\rho_1(\xi)$ $\mu$-a.e.  Hence,
 $\rho_1[\mathbb{E}_{\vX}(\gamma)](\omega)\in\varTheta(\omega)$ for every $\omega\in\vO$.
 It follows that
\[
X\cap\vD_{M_1\rho_1}(\omega)\subset{X}\cap\varTheta(\omega)\quad\mbox{for every }\omega\in\vO.
\]
If $\gamma\in\mcD_M^{mP}$, then $\mathbb{E}_{\vX}(\gamma)$ is strongly measurable and  Pettis integrable on $\vX$. Moreover, $\rho_1[\mathbb{E}_{\vX}(\gamma)]$  is strongly measurable with respect to $\vX$ and with almost all values in $X$. Since $\gamma\in{\mathcal {QS}_{\vG}}$ we have $\lg{x^*,\gamma}\rg\leq s(x^*,\vG)$ a.e. If $m(\bullet)=(P)\int_{\bullet}\mathbb{E}_{\vX}(\gamma)\,d\mu$ and $E\in\vX$, then
\begin{align*}
\lg{x^*,m(E)}\rg&=\int_E\lg{x^*,\mathbb{E}_{\vX}(\gamma)}\rg\,d\mu=\int_E\lg{x^*,\gamma}\rg\,d\mu\\
&\leq\int_E s(x^*,\vG)\,d\mu\stackrel{\eqref{E14}}{=}\int_Es(x^*,\mathbb{E}_{\vX}(\vG))\,d\mu=s(x^*,M_1(E))\,.
\end{align*}
Hence $m\in\mcS(M_1)$ and
$\xi:=\mathbb{E}_{\vX}(\gamma)\in\mcD_{M_1}^{mP}$. Consequently,
\[
X\cap\varTheta(\omega)\subset X\cap\vD_{M_1\rho_1}(\omega)\quad\mbox{for  every } \omega\in\vO\,.
 \]
\end{proof}
\begin{rem}\label{r5}\rm Assume that $X$ is separable, has RNP and $\vG:\vO\to{cb(X)}$ is measurable and Pettis integrable in $cb(X)$.  It is a consequence of \cite[Remark 3.12]{mu25} that if $\vD:=\mathbb{E}_{\vX}(\vG)$, then $\vD=\vD_{M_1\rho_1}^X=\vD_{M_1\rho_1}\cap{X}$ $\mu|_{\vX}$-a.e. Moreover, there exists a family $\{\gamma_n=\rho_1(\gamma_n):n\in\N\}$ of functions such that $\vD(\omega)=\ov{\{\gamma_n(\omega):n\in\N\}}$ $\mu|_{\vX}$-a.e. This fact will be applied in Remark \ref{r6} to obtain a new proof of \cite[Theorem 3.7]{EEE}.
\end{rem}
If the multimeasure $M_1$ above is only of $\sigma$-finite variation, then we obtain the following result:
\begin{thm}\label{t7}
Assume that $X$ has RNP, $\vG:\vO\to{c(X)}$ is  Pettis integrable in $cb(X)$ and $M:\vS\to{c(X)}$ is  its integral. Assume also that $\vX$ is a sub-$\sigma$-algebra of $\vS$ that is $\mu|_{\vX}$-complete, $\rho_1$ is a lifting on $\vX$ and the multimeasure $M$ restricted to $\vX$ is of $\sigma$-finite variation. Then there exists  a  measurable conditional expectation $\mathbb{E}_{\vX}(\vG):\vO\to{cb(X)}$ of $\vG$ with respect to $\vX$ satisfying for each open $U\subset{X}$ the inclusion
\begin{equation}
\{\omega\in\vO:\mathbb{E}_{\vX}(\vG)(\omega)\cap{U}\neq\emp\}\subset \rho_1\{\omega\in\vO:\mathbb{E}_{\vX}(\vG)(\omega)\cap{U}\neq\emp\}\,.
\end{equation}
More precisely,  the multifunctions defined by any of the following formulae:
\begin{align}
\mathbb{E}_{\vX}(\vG)(\omega)&=
\ov{ \left\{\rho_1^\#[\mathbb{E}_{\vX}(\gamma)](\omega)\colon \gamma\in\mcD_M^{mP}\right\}\cap{X}}\quad\mbox{if }\omega\in{Dom(\vG)}\label{E18}\\
\mathbb{E}_{\vX}(\vG)(\omega)&=
X\cap\ov{ \left\{\rho_1^\#[\mathbb{E}_{\vX}(\gamma)](\omega)\colon \gamma\in\mcD_M^{mP}\right\}}\quad\mbox{if } \omega\in{Dom(\vG)}\label{E19}
\end{align}
are conditional expectations of $\vG$ with respect to $\vX$ (please remember about \eqref{E51}) .
Moreover, if $x^*\in{X^*}$, then
\begin{equation}
s(x^*,\mathbb{E}_{\vX}(\vG)(\omega))\leq \rho^\#_1[s(x^*,\mathbb{E}_{\vX}(\vG))](\omega)\quad\mbox{for every }\omega\in{Dom(\vG)}.
\end{equation}
\end{thm}
\section{Martingales.}
We begin with a general convergence theorem for martingales valid without particular assumptions about the range Banach space. Several other results can be found in \cite{mu8}.
\begin{deff} \rm Let $(\vS_n)_n$ be a non-decreasing sequence of sub-$\sigma$-algebras of $\vS$. For each $n\in\N$ let $\vG_n:\vO\to{c(X)}$  be a scalarly $\vS_n$-measurable multifunction. The sequence $(\vG_n,\vS_n)_n$ is called a martingale if each $\vG_n$ is Pettis integrable in $c(X)$ and $\mathbb{E}_{\vS_n}(\vG_{n+1})\stackrel{s}{\equiv}\vG_n$. Equivalently, if $(P)\!\!\int_E\vG_{n+1}\,d\mu=(P)\!\!\int_E\vG_n\,d\mu$, for every $E\in\vS_n$.
\end{deff}
\begin{deff} \cite[p.852]{mu8} \rm A sequence $(\vG_n)_n$ of scalarly measurable multifunctions $\vG_n:\vO\to{c(X)}$ is called scalarly equi-Cauchy in measure if for every $\eta>0$
\[
\lim_{m,n\to\infty}\sup_{\|x^*\|\leq 1}\mu\{\omega\in\vO\colon |s(x^*,\vG_m(\omega))-s(x^*,\vG_n(\omega))|>\eta\}=0.
\]
Replacing $\vG_m$ by $\vG$ one obtains scalar equi-convergence of $(\vG_n)_n$ in measure to $\vG$.
\end{deff}
The subsequent result can be considered as a generalization of \cite[Theorem 4.6]{mu8}.
\begin{thm}\label{t2}
Let $\vG:\vO\to{c(X)}$ be a scalarly integrable multifunction and let $(\vS_n)_n$ be an increasing sequence of sub-$\sigma$-algebras of $\vS$ such that $\vS$ is  the $\mu$-completion of $\sigma\left(\bigcup_n\vS_n\right)$. Assume  that for each $n\in\N$ there exists a multifunction $\vG_n:\vO\to{c(X)}$ that is scalarly measurable with respect to $\vS_n$ and Pettis integrable in $cb(X)\;[cwk(X),ck(X)]$ on $\vS_n$. Assume also that for each $x^*\in{X^*}$ and $n\in\N$
\begin{equation}\label{E21}
\mathbb{E}_{\vS_n}(s(x^*,\vG))=s(x^*,\vG_n)\quad \mu-a.e.
\end{equation}
 or  $(\vG_n,\vS_n)_n$ is a martingale and
\begin{equation}\label{E32}
\lim_ns(x^*,\vG_n)=s(x^*,\vG)\quad \mbox{in }\mu-measure.
\end{equation}
If $\bigcup_n\mcZ_{\vG_n}$ is uniformly integrable,
then $\vG\in{\mathbb P}(\mu,cb(X))\;[cwk(X),ck(X)]$, $\mathbb{E}_{\vS_n}(\vG)\stackrel{s}{\approx}\vG_n$ for every $n\in\N$  and
\[
\lim_nd_H(M_{\vG}(E),M_{\vG_n}(E))=0\quad\mbox{for every }E\in\vS.
\]
Moreover,
\[
\lim_nd_P(\vG,\vG_n))=0
\]
if and only if the sequence $(\vG_n)_n$ is scalarly equi-Cauchy in measure.
\end{thm}
\begin{proof}
The sequence $(\vG_n,\vS_n)_n$ is a martingale and
due to the uniform integrability of $\bigcup_n\mcZ_{\vG_n}$, for each $x^*\in{X^*}$ the martingale $s((x^*,\vG_n),\vS_n)_n$ is bounded in $L_1(\mu)$ and so \eqref{E32} implies the convergence $\lim_n\int_{\vO}|s(x^*,\vG_n)-s(x^*,\vG)|\,d\mu=0$ and the equality $\mathbb{E}_{\vS_n}(s(x^*,\vG))=s(x^*,\vG_n)$ $\mu$-a.e. According to \cite[Proposition 3.2]{mu8} there exists an $h$-multimeasure $M:\vS\to{cb(X)}\;[cwk(X),ck(X)]$ with the property
\[
\lim_nd_H(M_{\vG_n}(E),M(E))=0\quad\mbox{for each }E\in\vS.
\]
Consequently, $\vG$ is Pettis integrable, $M=M_{\vG}$ and  $\mathbb{E}_{\vS_n}(\vG)\stackrel{s}{\approx}\vG_n$ for every $n\in\N$.

Assume that the sequence $(\vG_n)_n$ is scalarly equi-Cauchy in measure.
For each $x^*\in{X^*}$ and $n\in\N$ we have $\mathbb{E}_{\vS_n}(s(x^*,\vG))=s(x^*,\vG_n)$ $\mu$-a.e. and in $L_1(\mu)$.  Hence
\begin{equation*}
	\lim_n \int_{E} s(x^*, \vG_n)\, d\mu =  \int_{E}s (x^*, \vG)\, d\mu\quad\mbox{for each }E\in\vS.
\end{equation*}
Let us fix $E\in\vS$, $\varepsilon>0$ and $0<\eta<\ve$ . For each $n,m\in\N$ and $x^*\in B_{X^*}$ denote by $H_{mn,x^*}$ the set
\[H_{m,n,x^*}:= \{\omega \in\vO\colon |s(x^*,\vG_m(\omega))-s(x^*,\vG_n(\omega))|>\eta \}.\]
By the assumption the family $\bigcup_n\mcZ_{\vG_n}$ is uniformly integrable and so one may choose $\delta<\ve$ be such that if $\mu(F)<\delta$, then
\begin{equation}\label{E22}
\sup_n\sup_{\|x^*\|\leq 1}\int_F|s(x^*,\vG_n)|\,d\mu<\ve.
\end{equation}
The scalar equi-Cauchy condition in measure  yields existence of $k\in\N$ such that  for all $m,n\geq k$
\[
\sup_{\|x^*\|\leq 1} \mu(H_{m,n,x^*})<\delta.
\]
Let $\N\ni{n(x^*)}>k$ be such that $\sup_{n\geq{n(x^*)}}\int_{\vO}|s(x^*,\vG_n)-s(x^*,\vG)|\,d\mu<\ve$. Denote also -- for simplicity -- the set $H_{n,n(x^*),x^*}$ by $H_{nx^*}$.  Then, for all $n\geq{k}$ and a fixed $x^*\in{B_{X^*}}$
\begin{eqnarray*}
	&& \biggl|\int_E s(x^*,\vG_n)\,d\mu-\int_E s(x^*,\vG)\,d\mu  \biggr|   \\
&\leq& \biggl|\int_E s(x^*,\vG_n)\,d\mu-\int_E s(x^*,\vG_{n(x^*)})\,d\mu  \biggr|+
\biggl|\int_E s(x^*,\vG_{n(x^*)})\,d\mu-\int_E s(x^*,\vG)\,d\mu  \biggr|\\
	& \leq&  \biggl|\int_{E\cap H_{nx^*}}s(x^*,\vG_n)\,d\mu-\int_{E\cap H_{nx^*}}s(x^*,\vG_{n(x^*)})\,d\mu\biggr|+
	\\ &+&\biggl|\int_{E\cap H^c_{nx^*}}s(x^*,\vG_n)\,d\mu-\int_{E\cap H^c_{nx^*}}s(x^*,\vG_{n(x^*)})\,d\mu\biggr|+\ve.
\end{eqnarray*}

Observe that by (\ref{E22}), we have
\begin{eqnarray*}
	&& \biggl|\int_{E  \cap H_{n(x^*)}}s(x^*,\vG_n)\,d\mu-\int_{E\cap H_{n(x^*)}}s(x^*,\vG_m)\,d\mu\biggr| \leq \\
	&&   \leq \biggl|\int_{E \cap H_{n(x^*)}}s(x^*,\vG_n)\,d\mu \biggr|+
 \biggl|\int_{E\cap H_{n(x^*)} }s(x^*,\vG_m)\,d\mu\biggr| \leq 2\ve
\end{eqnarray*}
and then
\[
\biggl|\int_{E\cap{H^c_{n(x^*)}}} s(x^*,\vG_n)\, d\mu - \int_{E\cap{H^c_{n(x^*)}}} s(x^*,\vG_{n(x^*)}) \,d\mu \biggr| \leq   \eta \, \mu(E\cap{H^c_{n(x^*)}})<\ve.
\]
Consequently, if $n>k$, then
\[
\biggl|\int_E s(x^*,\vG_n)\,d\mu-\int_E s(x^*,\vG)\,d\mu  \biggr| <5\ve.
\]
Since $k$ is independent of $x^*$ and of $E\in\vS$, we have finally for $n>k$ and arbitrary $A\in\vS$
\begin{eqnarray*}
d_H(M_{\vG_n}(A), M(A))&=& \sup_{\|x^*\|\leq 1}|s(x^*, M_{\vG_n}(A)) -s(x^*,M(A))|  \\
&\leq &\sup_{\|x^*\|\leq 1}\biggl|\int_A s(x^*,\vG_n)\,d\mu-\int_A s(x^*,\vG)\,d\mu\biggr|\geq{5\ve}.
\end{eqnarray*}
An appeal to \eqref{E3} completes the proof.

If $\lim_nd_P(\vG,\vG_n)=0$, then obviously $(\vG_n)_n$ is scalarly equi-convergent in measure to $\vG$, so it is also scalarly equi-Cauchy in measure.
\end{proof}
\begin{rem}\rm Instead of  the uniform integrability of $\bigcup_n\mcZ_{\vG_n}$ one may assume that $\mcZ_{\vG}$ is uniformly integrable. Such an assumption yields the uniform integrability of $\bigcup_n\mcZ_{\vG_n}$ (see the proof of \cite[Proposition 3.3]{mu8}).

It is perhaps worth to recall here that contrary to Pettis integrable functions, if $\vG$ is Pettis integrable in $cb(X)$, then the set $\{s(x^*,\vG)\colon x^*\in{B_{X^*}}\}$ may be not uniformly integrable (see \cite[Remark 1.13]{mu4}).
It is however so in case of strongly Pettis integrable $\vG$ and in case of $X$ not containing any isomorphic copy of $c_0$. Indeed,
since $M_{\vG}$ is then an $h$-multimeasure (see \cite[Proposition 4.1]{cdpms2020} or \cite[Theorem 3.1]{dms}) its range $j\circ{M_{\vG}}$ in the Banach space $Y=\ell_{\infty}(B_{X^*})$ via the   R{\aa}dstr\"{o}m embedding $j:cb(X)\to\ell_{\infty}(B_{X^*})$  is a vector measure. Consequently, the family $\{\lg{y^*,j\circ{M_{\vG}}}\rg\colon y^*\in\ell^*_{\infty}(B_{X^*})\;\&\;\|y^*\|\leq 1\}$ is uniformly $\mu$-continuous and $\sup_{\|y^*\|\leq 1}\vert\lg{y^*,j\circ{M_{\vG}}}\rg\vert(\vO)<\infty$. It means in particular that the family $\mcZ_{\vG}$ is uniformly integrable. But then it follows from \cite[Proposition 3.3]{mu8} that also $\bigcup_n\mcZ_{\vG_n}$ is uniformly integrable and $\lim_nd_H(M_{\vG}(E),M_{\vG_n}(E))=0$ for each $E\in\vS$.
\end{rem}
 There is an amount of papers where martingale selections of a $cb(X)$-valued multimartingale are constructed in case of separable $X$ (cf. \cite{hess1},  \cite{CD}, \cite{EEE} ). Unfortunately, if $X$ is non-separable then  -- as a rule -- Castaing representation fails and in fact it is unknown whether even one martingale selection exists. When quasi-selections are taken into account, then an uncountable  Castaing type representations  of scalarly equivalent martingales become possible.
\begin{prop}\label{p8}
Assume that $X$ has RNP, $(\vG_n,\vS_n)_{n=1}^{\infty}$ is a martingale of multifunctions  $\vG_n\in{\mathbb P}(\mu,cb(X)),\;n\in\N$, and $\rho\in\vL(\mu)$. For each $n\in\N$ let $M_n$ be the Pettis integral of $\vG_n$.  Assume that every $\mu|_{\vS_n}$ is complete, $\vS$ is the $\mu$-completion of $\sigma(\bigcup_n\vS_n)$ and there exists a partition $\{F_k:k\in\N\}\subset\vS_1$ of $\vO$ such that
\begin{equation}\label{E23}
\sup_n|M_n|(F_k)<\infty\quad\mbox{for every }k\in\N\,.
\end{equation}
Then
\begin{itemize} 	
	 \item[\textbf{(\ref{p8}.i)}]	
there exists a measurable $cb(X)$-valued martingale $(\wt{\vG_n},\vS_n)_n$ that is scalarly equivalent to $(\vG_n,\vS_n)_n$;
\item[\textbf{(\ref{p8}.ii)}]
$(\wt{\vG_n},\vS_n)_n$ admits a family $\{\gamma_{nt}:t\in{T}\,,n\in\N\}$ of $X$-valued Pettis integrable functions such that for each $t\in{T}$ $(\gamma_{nt},\vS_n)_n$ is a martingale selection of $(\wt{\vG_n},\vS_n)_n$;
\item[\textbf{(\ref{p8}.iii)}]
$\wt{\vG_n}=\ov{conv}\{\gamma_{nt}:t\in{T}\}$ on $Dom(\vG_n)$, for each $n\in\N$;
\item[\textbf{(\ref{p8}.iv)}]
Each martingale selection of $(\wt{\vG_n},\vS_n)_n$ is. a.e. convergent. If  $\sup_n\|s(x^*,\vG_n)\|_{L_1}<\infty$ for every $x^*\in{X^*}$, then the limit function of every martingale selection is Pettis integrable.
\end{itemize}
\end{prop}
\begin{proof}
Since $\mu|_{\vS_n}$ is complete the restriction of $\rho$ to $\vS_n$ is a lifting on $\vS_n$.
It follows then from Proposition \ref{p7} that for each $n$ there exists a measurable with respect to $\vS_n$ multifunction $\wt{\vG_n}$ that is scalarly equivalent to $\vG_n$. That solves (\ref{p8}.i). Moreover,  $\wt\vG_n$ can be defined by any of the following equalities:
\begin{equation}\label{E24}
\wt\vG_n(\omega)=
\begin{cases}
\ov{\left\{\rho^\#(\gamma)(\omega): \gamma\in\mcD_{M_n}^{mP}\right\}\cap{X}}, &\text{if $\omega\in{Dom(\vG_n)}$;}\\
\\
X\cap\ov{\left\{\rho^\#(\gamma)(\omega):  \gamma\in\mcD_{M_n}^{mP}\right\}}, &\text{if $\omega\in{Dom(\vG_n)}$;}.
\end{cases}
\end{equation}
Let us fix $k\in\N$.  For each $m\in\N$  denote by $M_m^k$ the restriction of $M_m$ to the $\sigma$-algebra $\vX_m^k:=F_k\cap\vS_m$ on $F_k$.  If $m\geq{k}$ and $\gamma_{mk}\in\mcD_{M_m^k}^{mP}$, then there exists $\gamma_{{m+1},k}\in\mcD_{M_{m+1}^k}^P$ such that $\mathbb{E}_{\vX_m^k}(\gamma_{{m+1},k})=\gamma_{mk}$. That fact follows from RNP of $X$ and Proposition \ref{p3}, because then each $M_m^k$ is $\mu|_{\vX_m^k}$-rich in strong densities. $(\rho^\#(\gamma_{nk}),\vX_n^k)_{n\geq{k}}$ is a martingale selection of $(\wt{\vG_n}|_{F_k},\vX_n^k)_{n\geq{k}}$ bounded in $L_1(\mu|_{F_k},X)$. We extend the martingale for $m<k$ by setting $\gamma_{mk}:=\rho^\#[\mathbb{E}_{\vX_m^k}(\gamma_{kk})]$. If $\gamma_m:=\sum_k\gamma_{mk}\chi_{\rho(F_k)}$, then $\gamma_m$ is a scalarly integrable selection of $\wt\vG_m$. Since $c_0\nsubseteq{X}$ isomorphically, the strongly measurable function $\gamma_m$ is Pettis integrable (cf. \cite[Theorem II.3.7]{du}) and $(\gamma_m,\vX_m)_m$ defines a martingale. This proves (\ref{p8}.ii). The condition (\ref{p8}.iii) is a consequence of the construction of each $\wt\vG_m$ presented in Theorem \ref{t7}.

Since $X$ has RNP, each martingale  selection of $(\wt{\vG_n},\vS_n)_n$ is  convergent $\mu$-a.e. on $F_k$ to a strongly measurable function. Hence, the same is valid for the whole $\vO$. By the assumption each martingale  selection of $(\wt{\vG_n},\vS_n)_n$ is bounded in $L_1(\mu)$. Since $X$ does not contain any isomorphic copy of $c_0$ its limit function is Pettis integrable on $\vO$.
\end{proof}
\begin{rem}\label{r6}\rm  One can apply Remark \ref{r5}  to give a new proof of \cite[Theorem 3.7]{EEE}  without assuming RNP of $X^*$.
\end{rem}
\begin{rem}\rm The assumption \eqref{E23} can be weaken but for the price of a more complicated thesis. Namely, one could assume only that $F_k\in\vS_{n_k}$ with $n_k<n_{k+1}$ and $\sup_{n\geq{n_k}}|M_n|(F_k)<\infty$. Then, we apply the thesis of Proposition \ref{p8} for each $F_k$ separately and for each $k$ the set $T$ may be different.
\end{rem}
Contrary to Theorem \ref{t2}, assuming RNP of $X$ one gets a chance to prove existence of a limit multifunction invoking to properties of the martingale under consideration only. We need yet one more definition. We say that a sequence $(A_n)_n$ of elements of $cb(X)$ is $\tau_L$-converging to $A\in{cb(X)}$ if $s(x^*,A_n){\longrightarrow}{s(x^*,A)}$ for every $x^*\in{X^*}$ and $d(x,A_n)\longrightarrow{d(x,A)}$ for every $x\in{X}$. $d(x,B)$ is the distance of $x$ from $B$.

The subsequent result is a  generalization of \cite[Theorem 4.6]{EEE} in its existence part.
\begin{thm}\label{t8}
Assume that $X$ has RNP, $(\vG_n,\vS_n)_{n=1}^{\infty}$ is a martingale of multifunctions  $\vG_n:\vO\to{c(X)},\;n\in\N$, that are Pettis integrable in $cb(X)\,[cwk(X), ck(X)]$ and $\rho\in\vL(\mu)$. Assume that every $\mu|_{\vS_n}$ is complete, $\vS$ is the $\mu$-completion of $\sigma\left(\bigcup_n\vS_n\right)$ and there exists a partition $\{F_k:k\in\N\}\subset\vS_1$ of $\vO$ such that
\begin{itemize}
\item[\textbf{(\ref{t8}.i)}]
$\sup_n|M_{\vG_n}|(F_k)<\infty\quad\mbox{for every }k\in\N$\,;
\item[\textbf{(\ref{t8}.ii)}]
$\sup_n\|s(x^*,\vG_n)\|_{L_1}<\infty\quad\mbox{for every }x^*\in{X^*}$\,;
\end{itemize}
Then there exists a measurable multifunction $\vG\in{\mathbb P}(\mu,cb(X))\,[cwk(X), ck(X)]$ and a martingale $(\wt\vG_n,\vS_n)_n$ of measurable $cb(X)$-valued multifunctions such that
\begin{itemize}
\item[\textbf{(\ref{t8}.iii)}]
For each $k\in\N$\;\; $|M_{\vG}|(F_k)<\infty$;
\item[\textbf{(\ref{t8}.iv)}]
For each $x^*\in{X^*}$ \;$\lim_ns(x^*,\vG_n)=s(x^*,\vG)$ a.e. and $\lim_n\|s(x^*,\vG_n)-s(x^*,\vG)\|_{L_1}=0$;
\item[\textbf{(\ref{t8}.v)}]
For each $n\in\N$ $\wt\vG_n$ is a measurable version of $\mathbb{E}_{\vS_n}(\vG)$ and $\wt\vG_n\stackrel{s}{\approx}\vG_n$.
\end{itemize}
If moreover the set $\bigcup_n\mcZ_{\vG_n}$ is uniformly integrable, then
\begin{itemize}
\item[\textbf{(\ref{t8}.vi)}]
For each $E\in\vS$\; $\tau_L-lim_n(P)\!\!\int_E\vG_n\,d\mu=(P)\!\!\int_E\vG\,d\mu$;
\item[\textbf{(\ref{t8}.vii)}]
$\lim_nd_H(M_{\vG_n}(E),M_{\vG}(E))=0\quad for\; every \quad E\in\vS$.
\end{itemize}
\end{thm}
\begin{proof} In virtue of Proposition \ref{p8} there exists a martingale $(\gamma_n,\vS_n)_n$ such that each $\gamma_n$ is a Pettis integrable quasi-selection of $\vG_n$. Then $(\vG_n-\gamma_n,\vS_n)_n$ is a martingale of Pettis integrable multifunctions with the property $s(x^*,\vG_n-\gamma_n)\geq )$ $\mu$-a.e. for each $x^*\in{X^*}$ separately.  To simplify the proof we assume that we have already $s(x^*,\vG_n)\geq 0$ $\mu$-a.e.

\textbf{(\ref{t8}.iii)}-\textbf{(\ref{t8}.v)}.
Denote the algebra $\bigcup_n\vS_n$ by $\vS_0$. Without loss of generality we may also assume that $F_k=\rho(F_k),\;k\in\N$.
According to Proposition \ref{p8} the formula \eqref{E24} defines a martingale $(\wt\vG_n,\vS_n)_n$ that is scalarly equivalent to $(\vG_n,\vS_n)_n$ and each multifunction $\wt\vG_n$ is measurable with respect to $\vS_n$. In virtue of Proposition \ref{p7} the function $|\wt\vG_n|$ is measurable and so it  follows from \cite[Theorem 4.7]{mu4} that $|M_{\vG_n}|(E)=\int_E|\wt\vG_n|\,d\mu$, for every $E\in\vS_n$.

Let $k\in\N$ be fixed. It follows then from Fatou lemma and (\ref{t8}.i) that
\[
\int_{E\cap{F_k}}\liminf_n|\wt\vG_n|\,d\mu\leq\liminf_n\int_{E\cap{F_k}}|\wt\vG_n|\,d\mu<\infty\,.
\]
But that yields $\liminf_n|\wt\vG_n|\chi_{F_k}<\infty$ a.e. Let $\ve>0$ be arbitrary and let $\delta>0$ be such that $\mu(E\cap{F_k})<\delta$ yields $\int_{E\cap{F_k}}\liminf_n|\wt\vG_n|\,d\mu<\ve$. If $x^*\in{B_{X^*}}$, then $(s(x^*,\wt\vG_n))_n$ is a bounded martingale in $L_1(\mu|_{F_k})$  and so it is convergent in $L_1(\mu|_{F_k})$ and $\mu$-a.e. to a function in ${L_1(\mu|_{F_k})}$. Let
\[
M_{x^*}^k(E\cap{F_k}):=\int_{E\cap{F_k}}\lim_ns(x^*,\wt\vG_n))\,d\mu=\lim_n s(x^*,M_{\vG_n}(E\cap{F_k}))\quad\mbox{for }  E\in \vS_0\,.                \]
But $|s(x^*,\wt\vG_n)|\leq|\wt\vG_n|,\;n\in\N$, on $F_k$ and so $\lim_n|s(x^*,\wt\vG_n\chi_{F_k})|\leq\liminf_n|\wt\vG_n|\chi_{F_k}$ $\mu$-a.e.
Consequently, if $E\in\vS_0$ and $\mu(E\cap{F_k})<\delta$, then
\begin{align}\label{E25}
|M_{x^*}^k(E\cap{F_k})|&=\lim_n|s(x^*,M_{\vG_n}(E\cap{F_k}))|=\lim_n\left|\int_{E\cap{F_k}}s(x^*,\wt\vG_n)\,,d\mu\right|\\
&\leq\int_{E\cap{F_k}}\lim_n|s(x^*,\wt\vG_n)|\leq\int_{E\cap{F_k}}\liminf_n|\wt\vG_n|\,d\mu<\ve,\notag
\end{align}
for an arbitrary $x^*\in{B_{X^*}}$.
We have $M_{\vG_n}(E)=M_{\vG_{n+1}}(E)$ for every $n\in\N$ and $E\in\vS_n$ and so we can define $M^{\dagger}_k:\vS_0\cap{F_k}\to{cb(X)}$  by $M^{\dagger}_k(E):=M_{\vG_n}(E\cap{F_k})$ if $E\in\vS_n$. $M^{\dagger}_k$ is an additive set function of finite variation on  $F_k$.
Due to the initial assumption of positivity of the support functionals of every $\vG_n$, we have $0\in{M^{\dagger}_k(E)}$, for every $E\in\vS_0$.
One of the consequences of \eqref{E25} is absolute continuity of $M_{x^*}^k$ on the algebra $\vS_0\cap{F_k}$ in the $\ve-\delta$ sense. That means that each $M_{x^*}^k$ is $\sigma$-additive on the algebra $\vS_0\cap{F_k}$. But if $E\in\vS_0$, then
\[
M_{x^*}^k(E\cap{F_k})=\lim_n\int_{E\cap{F_k}} s(x^*,\vG_n)\,d\mu=\lim_ns(x^*,M_{\vG_n}(E\cap{F_k}))=s(x^*,M^{\dagger}_k(E\cap{F_k})).
 \]
 This gives countable additivity of $s(x^*,M^{\dagger}_k)$ on $\vS_0\cap{F_k}$. Moreover,  formula \eqref{E25} proves also that the $\ve-\delta$ absolute continuity of the family $\{s(x^*,M^{\dagger}_k)\colon \|x^*\|\leq 1\}$ is uniform on $B_{X^*}$. It follows from \cite[Proposition 2.5]{ka} that $M^{\dagger}_k$ has an extension to a multimeasure $M^k:\vS\cap{F_k}\to{cb(X)}$. $M^k$ is of finite variation on $F_k$ because $M^{\dagger}_k$ is of finite variation on $F_k$ and $\vS_0$ is $\mu$-dense in $\vS$.
$M^k$ is in fact countably additive in the Hausdorff metric of $cb(X)$. Moreover, $0\in{M^k(E\cap{F_k})}$ for every $E\in\vS$.

In virtue of Theorem \ref{t6}  $M^k$ can be represented as $M^k(E\cap{F_k})=(P)\!\!\int_{E\cap{F_k}}\vG^k\,d\mu$, where $\vG^k:F_k\to{cb(X)}$ is measurable with respect to $\vS\cap{F_k}$. But $M^k|_{\vS_n\cap{F_k}}=M_{\vG_n}|_{\vS_n\cap{F_k}}$, what yields (Theorem \ref{t7})
\begin{equation}\label{E26}
\mathbb{E}_{\vS_n}(\vG^k)\stackrel{s}{\approx}\vG_n|_{F_k}\stackrel{s}{\approx}\wt\vG_n|_{F_k}.
 \end{equation}
It follows that for each $x^*\in{X^*}$ we have also $s(x^*,\mathbb{E}_{\vS_n}(\vG^k))=s(x^*,\vG_n|_{F_k})=s(x^*,\wt\vG_n|_{F_k})$ $\mu$-a.e. But according to \eqref{E14}
 $s(x^*,\mathbb{E}_{\vS_n}(\vG^k))=\mathbb{E}_{\vS_n}s(x^*,\vG^k)$ $\mu$-a.e. and as $(\mathbb{E}_{\vS_n}s(x^*,\vG^k),\vS_n)_n$ is a martingale that is on $F_k$ $\mu$-a.e. and in $L_1(\mu|_{F_k})$ convergent to $s(x^*,\vG^k)$, we obtain the equality
\begin{equation}\label{E27}
\lim_ns(x^*,\vG_n|_{F_k})=\lim_ns(x^*,\wt\vG_n|_{F_k})=s(x^*,\vG^k)\quad\mu-\mbox{a.e. on }F_k\,.
\end{equation}
As  $0\in{M^k(E\cap{F_k})}$ for every $E\in\vS$, the sequence $\left(\bigoplus_{k=1}^m M^k(E\cap{F_k})\right)_m$ is non-decreasing. It follows from \cite[Proposition 1]{sz} that
\begin{align*}
\lim_m\sum_{k=1}^ms(x^*,M^k(E\cap{F_k}))=\lim_ms\left(x^*,\bigoplus_{k=1}^m M^k(E\cap{F_k})\right)=s\biggl(x^*,\ov{\bigcup_m\bigoplus_{k=1}^m M^k(E\cap{F_k})}\biggr)\,.
\end{align*}
We define now $M:\vS\to{c(X)}$ setting $M(E):=\ov{\bigcup_m\bigoplus_{k=1}^m M^k(E\cap{F_k})}$. Clearly $M|_{F_k}=M^k$ for every $k\in\N$.
One can easily check that
\begin{align}\label{e8}
s(x^*,M(E))&=\sup_m\sum_{k=1}^ms(x^*,M^k(E\cap{F_k}))\\
&=\lim_m\sum_{k=1}^ms(x^*,M(E\cap{F_k}))=\sum_{k=1}^{\infty}s(x^*,M(E\cap{F_k}))\notag
\end{align}
and the series is absolutely convergent. Assume now that $E=\bigcup_nE_n$ and the sets $E_n$ are pairwise disjoint. Then, applying \eqref{e8} and the countable additivity of $s(x^*,M)$ on each $F_k$, we have
\begin{align*}
s\biggl(x^*,M\biggl(\bigcup_nE_n\biggr)\biggr)&=\sum_{k=1}^{\infty}s\biggl(x^*,M\biggl(\bigcup_nE_n\cap{F_k}\biggr)\biggr)=\sum_{k=1}^{\infty}\sum_{n=1}^{\infty}s(x^*,M(E_n\cap{F_k}))\\ &=\sum_{n=1}^{\infty}\sum_{k=1}^{\infty}s(x^*,M(E_n\cap{F_k}))=\sum_{n=1}^{\infty}s(x^*,M(E_n)).
\end{align*}
Thus, $M$ is a $\mu$-continuous multimeasure of $\sigma$-finite variation.  $\vG:\vO\to{cb(X)}$ defined by $\vG|_{F_k}=\vG^k$ is an obvious candidate for the Pettis integrable density of $M$ with respect to $\mu$. Indeed, applying the absolute convergence of the series
$\sum_k\int_{E\cap{F_k}}|s(x^*,\vG^k)|\,d\mu$ following from \eqref{e8}, we have
\begin{align*}
s(x^*,M(E))&=\sum_ks(x^*,M(E\cap{F_k}))=\sum_k\int_{E\cap{F_k}}s(x^*,\vG^k)\,d\mu=\int_Es(x^*,\vG)\,d\mu.
\end{align*}
Thus, $\vG\in{\mathbb P}(\mu,cb(X))$ and $M=M_{\vG}$.
The scalar convergence required in (\ref{t8}.iv) is a consequence of \eqref{E27} and of the convergence of the scalar martingales in $L_1(\mu)$. (\ref{t8}.v) follows from \eqref{E26}.

\textbf{(\ref{t8}.vi)}-\textbf{(\ref{t8}.vii)}. \quad Assume now that  the set $\bigcup_n\mcZ_n$ is uniformly integrable. $M$ is rich in strong densities and so if $E\in\vS$, then $M(E)=\ov{\{(P)\!\int_E\gamma\,d\mu\colon \gamma\in\mcD^{mP}_M\}}$. If $\gamma$ is a Pettis integrable density of $m\in\mcS(M)$, then  $(\gamma_n:=\mathbb{E}_{\vS_n}(\gamma),\vS_n)_n$ is
 a martingale such that $\gamma_n$ is a quasi-selection of $\vG_n$. It follows from \cite{Ch} that the martingale is $\mu$-a.e. convergent to $\gamma$ and \cite[Corollary 5.3]{mu4} yields the convergence $\lim_n\|\gamma_n-\gamma\|_P=0$. Due to this convergence, if $x\in{X}$, then
\[
\limsup_nd\biggl(x,(P)\!\!\int_E\vG_n\,d\mu\biggr)\leq \lim_n\biggl\|x-(P)\!\!\int_E\gamma_n\,d\mu\biggr\|=\biggl\|x-(P)\!\!\int_E\gamma\,d\mu\biggr\|,
\]
what yields
\begin{equation}\label{E28}
\limsup_nd\biggl(x,(P)\!\!\int_E\vG_n\,d\mu\biggr)\leq d\biggl(x,(P)\!\!\int_E\vG\,d\mu\biggr).
\end{equation}
To finish the proof, we will apply \cite[Corollary 1.5.3]{Beer}: If $\emp\neq{B}\subset{X}$, then $d(x,B)=\sup\{\lg{x^*,x}\rg-s(x^*,B):x^*\in{B_{X^*}}\}$.
\begin{align}\label{E29}
d\biggl(x,(P)\!\!\int_E\vG\,d\mu\biggr)&=\sup\biggl\{\biggl|\lg{x^*,x}\rg-s\biggl(x^*,(P)\!\!\int_E\vG\,d\mu\biggr)\biggr|\colon x^*\in{B_{X^*}}\biggr\}\\
&=\sup\biggl\{\biggl|\lg{x^*,x}\rg-\lim_ns\biggl(x^*,(P)\!\!\int_E\vG_n\,d\mu\biggr)\biggr|\colon x^*\in{B_{X^*}}\biggr\}\notag\\
&=\sup\biggl\{\liminf_n\biggl|\lg{x^*,x}\rg-s\biggl(x^*,(P)\!\!\int_E\vG_n\,d\mu\biggr)\biggr|\colon x^*\in{B_{X^*}}\biggr\}\notag\\
&\leq\liminf_n\sup\biggl\{\biggl|\lg{x^*,x}\rg-s\biggl(x^*,(P)\!\!\int_E\vG_n\,d\mu\biggr)\biggr|\colon x^*\in{B_{X^*}}\biggr\}\notag\\
&=\liminf_n d\biggl(x,(P)\!\!\int_E\vG_n\,d\mu\biggr). \notag
\end{align}
The inequalities \eqref{E28} and \eqref{E29} give the required convergence in (\ref{t8}.vi).

The condition (\ref{t8}.vii) follows from  Theorem \ref{t2}.
\end{proof}

\end{document}